\documentclass[a4paper,12pt]{article}
\usepackage{amscd}
\usepackage{amsmath,amsfonts,amssymb,amscd}
\usepackage{indentfirst,graphicx,epsfig}
\usepackage{graphicx,psfrag}
\input{epsf}

\newtheorem{thm}{Theorem}

\newtheorem{lem}{Lemma}
\newtheorem{con}{Conjecture}

\def\qed{\hfill \nopagebreak\rule{5pt}{8pt}}
\def\pf{\noindent {\it Proof.} }

\title{\bf Randi\'c index, diameter and the average
distance\footnote{Supported by NSFC No.10831001, PCSIRT and the
``973" program. } }
\author{
\small Xueliang Li, Yongtang Shi\\
\small Center for Combinatorics and LPMC-TJKLC \\
\small Nankai University, Tianjin 300071, China \\
\small lxl@nankai.edu.cn; shi@cfc.nankai.edu.cn
\date{}}
\begin{document}
\maketitle
\begin{abstract}
The Randi\'c index of a graph $G$, denoted by $R(G)$, is defined as
the sum of $1/\sqrt{d(u)d(v)}$ over all edges $uv$ of $G$, where
$d(u)$ denotes the degree of a vertex $u$ in $G$. In this paper, we
partially solve two conjectures on the Randi\'c index $R(G)$ with
relations to the diameter $D(G)$ and the average distance $\mu(G)$
of a graph $G$. We prove that for any connected graph $G$ of order
$n$ with minimum degree $\delta(G)$, if $\delta(G)\geq 5$, then
$R(G)-D(G)\geq \sqrt 2-\frac{n+1} 2$; if $\delta(G)\geq n/5$ and
$n\geq 15$, $\frac{R(G)}{D(G)} \geq \frac{n-3+2\sqrt 2}{2n-2}$ and
$R(G)\geq \mu(G)$. Furthermore, for any arbitrary real number
$\varepsilon \ (0<\varepsilon<1)$, if $\delta(G)\geq \varepsilon n$,
then $\frac{R(G)}{D(G)} \geq \frac{n-3+2\sqrt 2}{2n-2}$ and
$R(G)\geq \mu(G)$ hold for sufficiently large $n$.\\
[2mm] Keywords: Randi\'c index; diameter; average distance; minimum degree\\
[2mm] AMS Subject Classification (2000): 05C12, 05C35, 92E10.
\end{abstract}

\section{Introduction}
The {\it Randi\'c index} $R(G)$ of a (molecular) graph $G$ was
introduced by the chemist Milan Randi\'c \cite{R} in 1975 as the sum
of $1/\sqrt{d(u)d(v)}$ over all edges $uv$ of $G$, where $d(u)$
denotes the degree of a vertex $u$ in $G$, i.e.,
$R(G)=\sum\limits_{uv\in E(G)}\frac 1 {\sqrt{d(u)d(v)}}$\,.
Recently, many results on the extremal theory of the Randi\'c index
have been reported (see \cite{LG}).

Given a connected, simple and undirected graph $G=(V,E)$ of order
$n$. The distance between two vertices $u$ and $v$ in $G$, denoted
by $d_G(u,v)$ (or $d(u,v)$ for short), is the length of a shortest
path connecting $u$ and $v$ in $G$. The $diameter$ $D(G)$ of $G$ is
the maximum distance $d(u,v)$ over all pairs of vertices $u$ and $v$
of $G$. The \textit{average distance} $\mu(G)$, an interesting
graph-theoretical invariant, is defined as the average value of the
distances between all pairs of vertices of $G$, i.e.,
$$\mu(G)=\frac{\sum_{u,v\in V}d(u,v)}{{n \choose 2}}.$$
For terminology and notations not given here, we refer to the book
of Bondy and Murty \cite{BM}.

There are many results on the relations between the Rand\'c index
and some other graph invariants, such as the minimum degree, the
chromatic number, the radius, and so on. In this paper, we will
consider the relations of the Randi\'c index with the diameter and
the average distance.

In \cite{AHZ1}, Aouchiche, Hansen and Zheng proposed the following
conjecture on the relation between the Randi\'c index and the
diameter.

\begin{con}[\cite{AHZ1}]\label{con1}
For any connected graph of order $n\geq 3$ with Randi\'c index
$R(G)$ and diameter $D(G)$,
$$R(G)-D(G)\geq \sqrt 2-\frac{n+1} 2 \qquad  and \qquad \frac{R(G)}{D(G)}
\geq \frac{n-3+2\sqrt 2}{2n-2},$$ with equalities if and only if
$G\cong P_n$.
\end{con}

In \cite{F}, Fajtlowicz proposed the following conjecture on the
relation between the Randi\'c index and the average distance.
\begin{con}[\cite{F}]\label{con2}
For all connected graphs $G$, $R(G)\geq \mu(G)$, where $\mu(G)$
denotes the average distance of $G$.
\end{con}

In the following, we will prove that for any connected graph $G$ of
order $n$ with minimum degree $\delta(G)$, if $\delta(G)\geq 5$,
then $R(G)-D(G)\geq \sqrt 2-\frac{n+1} 2$; if $\delta(G)\geq n/5$
and $n\geq 15$, $\frac{R(G)}{D(G)} \geq \frac{n-3+2\sqrt 2}{2n-2}$
and $R(G)\geq \mu(G)$. Furthermore, for any arbitrary real number
$\varepsilon \ (0<\varepsilon<1)$, if $\delta(G)\geq \varepsilon n$,
then $\frac{R(G)}{D(G)} \geq \frac{n-3+2\sqrt 2}{2n-2}$ and
$R(G)\geq \mu(G)$ hold for sufficiently large $n$.

\section{Main results}

At first, we recall some lemmas which will be used in the sequel.

\begin{lem}[Erd\"os et al. \cite{PPPT}]\label{lem1}
Let $G$ be a connected graph with $n$ vertices and minimum degree
$\delta(G)\geq 2$. Then $D(G)\leq \frac {3n} {\delta(G)+1}-1$.
\end{lem}

\begin{lem}[Kouider and Winkler \cite{KW}]\label{lem4}
If $G$ is a graph with $n$ vertices and minimum degree $\delta(G)$,
then the average distance satisfying $\mu(G)\leq \frac n
{\delta(G)+1}+2$.
\end{lem}

\begin{lem}[Li, Liu and Liu \cite{LL}]\label{lem2}
Let $G$ be a graph of order $n$ with minimum degree $\delta(G)=k$.
Then
$$ R(G)\geq\left\{
 \begin{array}{ll}
 \frac{k(k-1)}{2(n-1)}+\frac{k(n-k)}{\sqrt{k(n-1)}}  &\mbox{if $k\leq \frac n 2$}\\[3mm]
\frac{(n-p)(n-p-1)}{2(n-1)}+\frac{p(p+k-n)}{2k}+\frac{p(n-p)}{\sqrt{k(n-1)}}
&\mbox{if $k> \frac n 2$}
 \end{array}
 \right.
 $$
where $p$ is an integer given as follows:
$$ p=\left\{
 \begin{array}{ll}
 \frac n 2  &\mbox{if $n\equiv 0 \ (mod~ 4)$}\\[1mm]
\left\lfloor\frac n 2 \right\rfloor \mbox{or} \left\lceil\frac n 2
\right\rceil &\mbox{if $n\equiv 1 \ (mod~ 4)$ and $k$ is even}\\[1mm]
\left\lfloor\frac n 2 \right\rfloor  &\mbox{if $n\equiv 1 \ (mod~ 4)$ and $k$ is odd}\\[1mm]
\frac {n-2} 2~  or ~\frac {n+2} 2 &\mbox{if $n\equiv 2 \ (mod~ 4)$ and $k$ is even}\\[1mm]
\frac n 2  &\mbox{if $n\equiv 2 \ (mod~ 4)$ and $k$ is odd}\\[1mm]
\left\lfloor\frac n 2 \right\rfloor \mbox{or} \left\lceil\frac n 2
\right\rceil &\mbox{if $n\equiv 3 \ (mod~ 4)$ and $k$ is even}\\[1mm]
\left\lceil\frac n 2
\right\rceil &\mbox{if $n\equiv 3 \ (mod~ 4)$ and $k$ is odd}.\\[1mm]
 \end{array}
 \right.
 $$
\end{lem}
It is easy to see from Lemma \ref{lem2} that $p$ is among the
numbers $\frac {n-2} 2$, $\frac {n-1} 2$, $\frac {n} 2$, $\frac
{n+1} 2$ and $\frac {n+2} 2$.

\begin{lem}\label{lem31}
Denote by $g(n,k)=\frac{2k-1}{2(n-1)}+\frac{n-3k}{2\sqrt{k(n-1)}}$.
Then for $1\leq k\leq n/2$, $g(n,k)\geq 0$.
\end{lem}
\pf If $n\geq 3k$, we can directly obtain that $g(n,k)>0$. Now we
assume that $2k\leq n<3k$. Then $$g(n,k)=\frac 1
{2\sqrt{n-1}}\left(\frac
{2k-1}{\sqrt{n-1}}-\frac{3k-n}{\sqrt{k}}\right)=\frac 1
{2\sqrt{n-1}}\left(\frac {2k-1}{\sqrt{n-1}}-3\sqrt{k}+\frac n
{\sqrt{k}}\right).$$ Since
\begin{align*}
\frac{\partial g(n,k)}{\partial k}&=\frac 1 {2\sqrt{n-1}}\left(\frac
2 {\sqrt{n-1}}-\frac 3 {2\sqrt{k}}-\frac n {2k\sqrt{k}}\right)\\
&<\frac 1 {2\sqrt{n-1}}\left(\frac 2 {\sqrt{n-1}}-\frac 3
{2\sqrt{k}}-\frac {2k} {2k\sqrt{k}}\right)<0,
\end{align*}
for $n\geq 2$, we have
$$g(n,k)>g(n,\frac n 2)=\sqrt{n-1}-\sqrt{\frac n 2}\geq 0.$$
Therefore, the lemma follows.\qed

\begin{thm}
For any connected graph $G$ of order $n$ with minimum degree
$\delta(G)$.

(1) \ If $\delta(G)\geq 5$, then $R(G)-D(G)\geq \sqrt 2-\frac{n+1}
2$;

(2) \ If $\delta(G)\geq n/5$ and $n\geq 15$, then $\frac{R(G)}{D(G)}
\geq \frac{n-3+2\sqrt 2}{2n-2}$. Furthermore, for any arbitrary real
number $\varepsilon \ (0<\varepsilon<1)$, if $\delta(G)\geq
\varepsilon n$, then $\frac{R(G)}{D(G)} \geq \frac{n-3+2\sqrt
2}{2n-2}$ holds for sufficiently large $n$.

(3) \ If $\delta(G)\geq n/5$ and $n\geq 15$, then $R(G)\geq \mu(G)$.
Furthermore, for any arbitrary real number $\varepsilon \
(0<\varepsilon<1)$, if $\delta(G)\geq \varepsilon n$, then $R(G)\geq
\mu(G)$ holds for sufficiently large $n$.
\end{thm}
\pf Let $G$ be a connected graph of order $n$ with minimum degree
$\delta(G)=k$. By Lemma \ref{lem1}, we have $D(G)\leq \frac {3n}
{k+1}-1$.

(1) \ Suppose $k\geq 5$, we will show $R(G)-D(G)\geq \sqrt
2-\frac{n+1} 2$. We consider the following two cases:

{\bf Case 1.} $k\leq \frac n 2$.

By Lemma \ref{lem2}, we only need to consider the following
inequality,
$$\frac{k(k-1)}{2(n-1)}+\frac{k(n-k)}{\sqrt{k(n-1)}}\geq \frac {3n}
{k+1}-1+\sqrt 2-\frac{n+1} 2.$$
Let
$$f(n,k)=\frac{k(k-1)}{2(n-1)}+\frac{k(n-k)}{\sqrt{k(n-1)}}-\frac
{3n} {k+1}-\sqrt{2}+\frac {n+3} 2.$$ Then by Lemma \ref{lem31},
$\frac{\partial f(n,k)}{\partial
k}>\frac{2k-1}{2(n-1)}+\frac{n-3k}{2\sqrt{k(n-1)}}>0$. If $k\geq 5$,
we have $f(n,k)\geq f(n, 5)=\frac {10} {n-1}
+\frac{5(n-5)}{\sqrt{5(n-1)}}-\sqrt{2}+\frac 3 2>0$.

{\bf Case 2.}  $\frac n 2< k\leq n-1 $.

Let
$q(n,p)=\frac{(n-p)(n-p-1)}{2(n-1)}+\frac{p(p+k-n)}{2k}+\frac{p(n-p)}{\sqrt{k(n-1)}}$.
In the following, we will show that for every $p\in\{\frac {n-2}
2,\frac {n-1} 2,\frac {n} 2,\frac {n+1} 2,\frac {n+2} 2\}$,
$$q(n,p)\geq \frac {3n}
{k+1}-1+\sqrt 2-\frac{n+1} 2.$$
In fact, if $p=\frac{n-2}{2}$,
denote by
\begin{align*}
h(n,k)=&~q(n,\frac{n-2}{2})-\frac
{3n} {k+1}-\sqrt{2}+\frac {n+3} 2\\
=&~\frac{n(n+2)}{8(n-1)}+\frac{(n-2)(2k-(n+2))}{8k}+\frac{n^2-4}{4\sqrt{k(n-1)}}\\
&~-\frac {3n} {k+1}-\sqrt{2}+\frac {n+3} 2.
\end{align*}
Notice that
$$
\frac {\partial h(n,k)} {\partial k}
=\frac{n^2-4}{8k^2}-\frac{n^2-4}{8k\sqrt{k(n-1)}}+\frac{3n)}{(k+1)^2}>0,
$$
since $8k^2\leq 8k\sqrt{k(n-1)}$, i.e., $\sqrt{k}\leq \sqrt{n-1}$.
Thus, we have
\begin{align*}
h(n,k)> h(n,\frac n
2)&=\frac{n(n+2)}{8(n-1)}-\frac{n-2}{2n}+\frac{n^2-4}{2\sqrt{2n(n-1)}}-\frac
{6n} {n+2}-\sqrt{2}+\frac {n+3} 2\\
&>\frac n 8+\frac{n-2} 8-\frac n 4+\frac{n^2-4}
n-6-\sqrt{2}+\frac{n+3} 2\\
&=\frac{2n+5} 4+\frac{n^2-4}{2\sqrt{2}}-6-\sqrt{2}.
\end{align*}
By some calculations, we have that $\frac{2n+5}
4+\frac{n^2-4}{2\sqrt{2}}-6-\sqrt{2}>0$ for $n\geq 8$. For $4\leq
n\leq 7$, it is easy to verify $h(n,\frac n 2)>0$.

 In a similar way, we can verify the inequality for each of the
cases for $p=\frac {n-1} 2$, $\frac {n} 2$, $\frac {n+1} 2$ or
$\frac {n+2} 2$. The details are omitted.

(2) \ Similarly, we consider the following two cases:

{\bf Case 1.} $k\leq \frac n 2$.

By Lemma \ref{lem2}, we only need to consider the following
inequality,
$$\frac{k(k-1)}{2(n-1)}+\frac{k(n-k)}{\sqrt{k(n-1)}}\geq \left(\frac {3n}
{k+1}-1\right)\frac{n-3+2\sqrt 2}{2n-2}.$$ Let
$$f(n,k)=\frac{k(k-1)}{2(n-1)}+\frac{k(n-k)}{\sqrt{k(n-1)}}-\left(\frac {3n}
{k+1}-1\right)\frac{n-3+2\sqrt 2}{2n-2}.$$ Then by Lemma
\ref{lem31}, $\frac{\partial f(n,k)}{\partial
k}>\frac{2k-1}{2(n-1)}+\frac{n-3k}{2\sqrt{k(n-1)}}>0$.

If $k\geq \frac n 5$ and $n\geq 15$, we have $f(n,k)\geq f(n, \frac
n 5)=\frac {n(n-5)} {50(n-1)}+\frac{4n^2}{5\sqrt{5n(n-1)}}-\frac
{(14n-5)(n-3+2\sqrt{2})}{2(n-1)(n+5)}>0$. Actually, for any
arbitrary positive number $\varepsilon \ (0< \varepsilon<1)$, if
$k\geq \varepsilon n$, then $f(n,k)>f(n,\varepsilon
n)>\frac{\varepsilon(1-\varepsilon)n^2}{\sqrt{\varepsilon
n(n-1)}}-\left(\frac {3n} {\varepsilon n+1}-1\right)\frac{n-3+2\sqrt
2}{2n-2}>0$ for sufficiently large $n$.

{\bf Case 2.}  $\frac n 2< k\leq n-1 $.

Let
$q(n,p)=\frac{(n-p)(n-p-1)}{2(n-1)}+\frac{p(p+k-n)}{2k}+\frac{p(n-p)}{\sqrt{k(n-1)}}$.
In the following, we will show that for every $p\in\{\frac {n-2}
2,\frac {n-1} 2,\frac {n} 2,\frac {n+1} 2,\frac {n+2} 2\}$,
$$q(n,p)\geq \left(\frac {3n}
{k+1}-1\right)\frac{n-3+2\sqrt 2}{2n-2}.$$ In fact, if
$p=\frac{n-2}{2}$, denote by
\begin{align*}
h(n,k)=&~q(n,\frac{n-2}{2})-\left(\frac {3n}
{k+1}-1\right)\frac{n-3+2\sqrt 2}{2n-2}\\
=&~\frac{n(n+2)}{8(n-1)}+\frac{(n-2)(2k-(n+2))}{8k}+\frac{n^2-4}{4\sqrt{k(n-1)}}\\
&~-\left(\frac {3n} {k+1}-1\right)\frac{n-3+2\sqrt 2}{2n-2}.
\end{align*}
Notice that
$$
\frac {\partial h(n,k)} {\partial k}
>\frac{n^2-4}{8k^2}-\frac{n^2-4}{8k\sqrt{k(n-1)}}>0,
$$
since $8k^2\leq 8k\sqrt{k(n-1)}$, i.e., $\sqrt{k}\leq \sqrt{n-1}$.
Thus, we have
\begin{align*}
h(n,k)> h(n,\frac n
2)&=\frac{n(n+2)}{8(n-1)}-\frac{n-2}{2n}+\frac{n^2-4}{2\sqrt{2n(n-1)}}\\
&~~~~-\left(\frac {6n}
{n+2}-1\right)\frac{n-3+2\sqrt 2}{2n-2}\\
&>\frac n 8+\frac{n-2} 8-\frac n 4+\frac{n^2-4}
n-\frac {(6-1)n} {2n-2}\\
&=\frac{n^2-4}{n}-\frac{11n-1}{4(n-1)}.
\end{align*}
By some calculations, we have that
$\frac{n^2-4}{n}-\frac{11n-1}{4(n-1)}>0$ for $n\geq 5$.

 In a similar way, we can verify the inequality for each of the
cases for $p=\frac {n-1} 2$, $\frac {n} 2$, $\frac {n+1} 2$ or
$\frac {n+2} 2$. The details are omitted.

By the method similar to (2), we can obtain the result of (3).

The proof is now complete.\qed

\end{document}